\documentclass[11pt]{article}
\usepackage{amstex,amssymb}

\newcommand{\qed}{\rule{3mm}{3mm}}

\newcommand{\ed}{{\bf 1}}

\newcommand{\bL}{{\bf L}}

\newcommand{\bS}{{\bf S}}

\newcommand{\bbL}{{\mathbb L}}
\newcommand{\bbS}{{\mathbb S}}

\newcommand{\cL}{{\cal L}}

\newcommand{\cX}{{\cal X}}

\newcommand{\gog}{{\mathfrak g}}

\newtheorem{theorem}{Theorem}[section]
\newtheorem{proposition}[theorem]{Proposition}
\newtheorem{lemma}[theorem]{Lemma}

\begin{document}
%
\begin{center}
{\large\bf 
Discrete Lagrangian reduction, discrete Euler--Poincar\'e equations,\\
and semidirect products}\end{center}
\vspace{0.5cm}

\begin{center}
{\sc Alexander I.\,Bobenko}\footnote{E--mail: 
{\tt bobenko} {\makeatother @@ \makeatletter}{\tt math.tu-berlin.de }} and 
{\sc Yuri B.\,Suris}\footnote{E--mail: 
{\tt suris} {\makeatother @@ \makeatletter}{\tt sfb288.math.tu-berlin.de }}
\end{center}
\begin{center}
Fachbereich Mathematik, Technische Universit\"at Berlin, \\
Str. 17 Juni 136, 10623 Berlin, Germany
\end{center}
\vspace{0.5cm}

\begin{abstract}
A discrete version of Lagrangian reduction is developed in the context of
discrete time Lagrangian systems on $G\times G$, where $G$ is a Lie group.
We consider the case when the Lagrange function is invariant with respect to the
action of an isotropy subgroup of a fixed element in the representation
space of $G$. In this context the reduction of the discrete Euler--Lagrange
equations is shown to lead to the so called discrete Euler--Poincar\'e
equations. A constrained variational principle is derived.
The Legendre transformation of the discrete Euler--Poincar\'e equations
leads to discrete Hamiltonian (Lie--Poisson) systems on a dual space 
to a semiproduct Lie algebra.
\end{abstract}

\newpage

\section{Introduction}
Dynamical systems with symmetry play an important role in mathematical
modelling of a vast variety of physical and mechanical processes. A Hamiltonian
approach to such systems is nowadays a well--established theory [A], [AM], [CB],
[MR], [MRW], [RSTS]. More recently, also a variational (Lagrangian) 
description of systems with symmetries also attracted much attention [MS],
[HMR], [CHMR]. In particular, in the last two papers the corresponding theory
was developed for Lagrangian systems on Lie groups, i.e. for Lagrangians
defined on tangent bundles $TG$ of Lie groups. A symmetry of the Lagrangian with
respect to a subgroup action leads to a reduced system on a semidirect
product described by the so called Euler--Poincar\'e equation.

In the present paper we develop a discrete analog of this theory, i.e. for
Lagrangians defined on $G\times G$. We introduce the corresponding reduced
systems and derive the discrete Euler--Poincar\'e equations. We establish
symplectic properties of the corresponding discrete dynamical systems. The
continuous time theory may be considered as a limiting case of the discrete 
time one. The important particular case, when the representation of $G$ 
participating in the general theory is chosen to be the adjoint representation,
is developed in [BS].

\setcounter{equation}{0}
\section{Lagrangian mechanics on $TG$ and on $G\times G$}

Recall that a continuos time Lagrangian system is defined by a smooth function
$\bL(g,\dot{g})\,:\,TG\mapsto{\mathbb R}$ on the tangent bundle of a smooth
manifold $G$. The function $\bL$ is called the {\it Lagrange function}. We
will be dealing here only with the case when $G$ carries an additional
structure of a {\it Lie grioup}. For an arbitrary function 
$g(t)\,:\,[t_0,t_1]\mapsto G$ one can consider the {\it action functional}
\begin{equation}\label{action}
\bS=\int_{t_0}^{t_1}\bL(g(t),\dot{g}(t))dt\;.
\end{equation} 
A standard argument shows that the functions $g(t)$ yielding extrema of this
functional (in the class of variations preserving $g(t_0)$ and $g(t_1)$), 
satisfy with necessity the {\it Euler--Lagrange equations}. In local 
coordinates $\{g^i\}$ on $G$ they read:
\begin{equation}\label{EL gen}
\frac{d}{dt}\left(\frac{\partial\bL}{\partial\dot{g}^i}\right)=
\frac{\partial\bL}{\partial g^i}\;.
\end{equation}
The action functional $S$ is independent of the choice of local coordinates,
and thus the Euler--Lagrange equations are actually coordinate independent as
well. For a coordinate--free description in the language of differential
geometry, see \cite{A}, \cite{MR}.

Introducing the quantities 
\begin{equation}\label{Pi}
\Pi=\nabla_{\dot{g}}\bL\in T_g^* G\;,
\end{equation}
one defines the {\it Legendre transformation}:
\begin{equation}\label{Legendre gen}
(g,\dot{g})\in TG\mapsto (g,\Pi)\in T^*G\;.
\end{equation}
If it is invertible, i.e. if $\dot{g}$ can be expressed through $(g,\Pi)$, 
then the the Legendre transformation of the Euler--Lagrange equations 
(\ref{EL gen}) yield  a {\it Hamiltonian system} on $T^*G$ with respect to the 
standard symplectic structure on $T^*G$ and with the Hamilton function
\begin{equation}\label{Ham gen}
H(g,\Pi)=\langle \Pi,\dot{g}\rangle-\bL(g,\dot{g})\;,
\end{equation}
(where, of course, $\dot{g}$ has to be expressed through $(g,\Pi)$). 

We now turn to the discrete time analog of these constructions,
introduced in \cite{V}, \cite{MV}. Our presentation is an adaptation of
the Moser--Veselov construction for the case when the basic manifold is a 
Lie group. Almost all constructions and results of the continuous time
Lagrangian mechanics have their discrete time analogs. The only exception 
is the existence of the ``energy'' integral (\ref{Ham gen}).

Let $\bbL(g,\hat{g})\,:\,G\times G$ be a smooth function, called the (discrete
time) {\it Lagrange function}. For an arbitrary sequence $\{g_k\in G,\;k=k_0,
k_0+1,\ldots,k_1\}$ one can consider the {\it action functional}
\begin{equation}\label{dS}
\bbS=\sum_{k=k_0}^{k_1-1}\bbL(g_k,g_{k+1})\;.
\end{equation}
Obviously, the sequences $\{g_k\}$ delivering extrema of this functional
(in the class of variations preserving $g_{k_0}$ and $g_{k_1}$), satisfy
with necessity the {\it discrete Euler--Lagrange equations}:
\footnote{For the notations from the Lie groups theory used in this and 
subsequent sections see Appendix B. In particular, for an arbitrary smooth 
function $f:G\mapsto{\mathbb R}$ its right Lie derivative $d\,'f$ and left Lie 
derivative $df$ are functions from $G$ into $\gog^*$ defined via the formulas
\[
\langle df(g),\eta\rangle=\left.\frac{d}{d\epsilon}\,f(e^{\epsilon\eta}g)
\right|_{\epsilon=0}\;,\qquad 
\langle d\,'f(g),\eta\rangle=\left.\frac{d}{d\epsilon}\,f(ge^{\epsilon\eta})
\right|_{\epsilon=0}\;,\qquad \forall \eta\in\gog\;,
\]
and the gradient $\nabla f(g)\in T^*_gG$ is defined as
\[
\nabla f(g)=R^*_{g^{-1}}\,df(g)=L^*_{g^{-1}}\,d\,'f(g).
\]
}
\begin{equation}\label{dEL}
\nabla_1\bbL(g_k,g_{k+1})+\nabla_2\bbL(g_{k-1},g_k)=0\;.
\end{equation}
Here $\nabla_1\bbL(g,\hat{g})$ ($\nabla_2\bbL(g,\hat{g})$) denotes the 
gradient of $\bbL(g,\hat{g})$ with respect to the first argument $g$ 
(resp. the second argument $\hat{g}$). So, in our case, when $G$ is a Lie 
group and not just a general smooth manifold, the equation (\ref{dEL}) is 
written in a coordinate free form, using the intrinsic notions of the Lie 
theory. As pointed out above, an invariant formulation of the Euler--Lagrange 
equations in the continuous time case is more sophisticated. This seems 
to underline the fundamental character of discrete Euler--Lagrange equations. 

The equation (\ref{dEL}) is an implicit equation for $g_{k+1}$. In 
general, it has more than one solution, and therefore defines a correspondence
(multi--valued map)  $(g_{k-1},g_k)\mapsto(g_k,g_{k+1})$. To discuss symplectic
properties of this correspondence, one defines:
\begin{equation}\label{dPi}
\Pi_k=\nabla_2\bbL(g_{k-1},g_k)\in T^*_{g_k}G\;.
\end{equation} 
Then (\ref{dEL}) may be rewritten as the following system:
\begin{equation}\label{dEL syst}
\left\{\begin{array}{l}
\Pi_k=-\nabla_1\bbL(g_k,g_{k+1}) \\ \\ \Pi_{k+1}=\nabla_2\bbL(g_k,g_{k+1}) 
\end{array}\right.
\end{equation}
This system defines a (multivalued) map $(g_k,\Pi_k)\mapsto(g_{k+1},\Pi_{k+1})$
of $T^*G$ into itself. More precisely, the first equation in (\ref{dEL syst})
is an implicit equation for $g_{k+1}$, while the second one allows for the
explicit and unique calculation of $\Pi_{k+1}$, knowing $g_k$ and $g_{k+1}$. As 
demonstrated in \cite{V}, \cite{MV}, this map $T^*G\mapsto T^*G$ is symplectic
with respect to the standard symplectic structure on $T^*G$.

\setcounter{equation}{0}
\section{Left symmetry reduction}

We want to consider the Lagrangian reduction procedure, in the case when
the Lagrange function is symmetric with respect to the action of a certain
subgroup of $G$ (precise formulations will follow). It turns out to be
convenient to perform this reduction in two steps. The first one of them
is not related to any symmetry, and is quite general.
 
\subsection{Left trivialization}
The tangent bundle $TG$ does not appear in the discrete time context at all. 
On the contrary, the cotangent bundle $T^*G$ still plays an important role in 
the discrete time theory, as the phase space with the canonical invariant 
symplectic structure. When working with the cotangent bundle of the Lie group, it is 
convenient to trivialize it, translating all covectors to the group unit by
left or right multiplication. This subsection is devoted to the constructions 
related to the left trivialization of the cotangent bundle $T^*G$:
\begin{equation}\label{d left triv *}
(g_k,M_k)\in G\times\gog^* \;\mapsto\; (g_k,\Pi_k)\in T^*G\;,
\end{equation}
where 
\begin{equation}\label{d M}
\Pi_k=L_{g_k^{-1}}^* M_k \quad\Leftrightarrow\quad M_k=L_{g_k}^*\Pi_k\;.
\end{equation}
Consider also the map
\begin{equation}\label{d left triv}
(g_k,W_k)\in G\times G \;\mapsto\; (g_k,g_{k+1})\in G\times G\;,
\end{equation}
where
\begin{equation}\label{W}
g_{k+1}=g_kW_k \quad\Leftrightarrow\quad W_k=g_k^{-1}g_{k+1}\;.
\end{equation}
In the continuous limit the elements $W_k$ lie in a  neighborhood of the 
group unit $e$.

Denote the pull--back of the Lagrange function under (\ref{d left triv}) through
\begin{equation}\label{dLagr left}
\bbL^{(l)}(g_k,W_k)=\bbL(g_k,g_{k+1})\;.
\end{equation}
We want to find difference equations satisfied by the sequences
$\{(g_k,W_k)\,,k=k_0,\ldots,k_1-1\}$ delivering extrema of the action 
functional
\begin{equation}\label{funct left triv}
\bbS^{(l)}=\sum_{k_0}^{k_1-1}\bbL^{(l)}(g_k,W_k)\;,
\end{equation}
and satisfying $W_k=g_k^{-1}g_{k+1}$.
Admissible variations of $\{(g_k, W_k)\}$ are those preserving the values 
of $g_{k_0}$ and $g_{k_1}=g_{k_1-1}W_{k_1-1}$. A more explicit description
of admissible variations is given by the following statement.
\begin{lemma}
The original variational problem for the functional $\bbS$ {\rm(\ref{dS})}
is equivalent to finding extremals of the functional $\bbS^{(l)}$
{\rm(\ref{funct left triv})} with admissible variations $\{(\tilde{g}_k,
\tilde{W}_k)\}$ of the sequence $\{(g_k,W_k)\}$ of the form
\begin{equation}\label{variation left}
\tilde{g}_k=g_k e^{\eta_k}\;, \qquad \tilde{W}_k=W_k e^{\eta_{k+1}-
{\rm Ad}\, W_k^{-1}\cdot \eta_k}\;,
\end{equation}
where $\{\eta_k\}_{k=k_0}^{k_1}$ is an arbitrary sequence of elements of the
Lie algebra $\gog$ with $\eta_{k_0}=\eta_{k_1}=0$.
\end{lemma}
{\bf Proof.} Obviously, the first formula in (\ref{variation left}) gives
a generic variation of $g_k$. We have: 
\[
\tilde{W}_k=\tilde{g}_k^{-1}\tilde{g}_{k+1}=
e^{-\eta_k}W_ke^{\eta_{k+1}}=W_ke^{-{\rm Ad}\, W_k^{-1}\cdot
\eta_k}e^{\eta_{k+1}}\;.
\]
Supposing now that all $\eta_k$ are small (of order $\epsilon$), we find, in
the first order in $\epsilon$, the second formula in (\ref{variation left}).
(Clearly, for variational purposes only the first order terms are essntial.)
\qed
\vspace{2mm}

\begin{proposition}
The difference equations for extremals of the functional $\bbS^{(l)}$ read:
\begin{equation}\label{dEL left Ham}
\left\{\begin{array}{l}
{\rm Ad}^*\,W_k^{-1}\cdot M_{k+1}=M_k+d\,'_{\!g}\bbL^{(l)}(g_k,W_k)\;,\\ \\
g_{k+1}=g_kW_k\;,\end{array}\right.
\end{equation}
where
\begin{equation}\label{d M thru L}
M_k=d\,'_{\!W}\bbL^{(l)}(g_{k-1},W_{k-1})\in\gog^*\;.
\end{equation}
If the ``Legendre transformation''
\begin{equation}\label{dLegendre left}
(g_{k-1},W_{k-1})\in G\times G\;\mapsto\; (g_k,M_k)\in G\times\gog^*\;,
\end{equation}
where $g_k=g_{k-1}W_{k-1}$, is invertible, then {\rm(\ref{dEL left Ham})} 
defines a map $(g_k,M_k)\mapsto(g_{k+1},M_{k+1})$ which is symplectic 
with respect to the following Poisson bracket on $G\times\gog^*$:
\begin{equation}\label{PB left}
\{f_1,f_2\}=-\langle d\,'_{\!g}f_1,\nabla_M f_2\rangle+
\langle d\,'_{\!g}f_2,\nabla_M f_1\rangle+
\langle M,[\nabla_M f_1,\nabla_M f_2]\,\rangle\;.
\end{equation}
\end{proposition}
{\bf The first proof.} The simplest way to derive (\ref{dEL left Ham}) is to 
pull back the equations (\ref{dEL}) under the map (\ref{d left triv}).
To do this, first rewrite (\ref{dEL}) as
\begin{equation}\label{dEL left aux1}
d\,'_{\!1}\bbL(g_k,g_{k+1})+d\,'_{\!2}\bbL(g_{k-1},g_k)=0\;.
\end{equation}
We have to express these Lie derivatives in terms of $(g,W)$. 
The answer is this:
\begin{equation}\label{dEL left aux2}
d\,'_{\!2}\bbL(g_{k-1},g_k)=d\,'_{\!W}\bbL^{(l)}(g_{k-1},W_{k-1})\;,
\end{equation}
\begin{equation}\label{dEL left aux3}
d\,'_{\!1}\bbL(g_k,g_{k+1})=d\,'_{\!g}\bbL^{(l)}(g_k,W_k)-
d_W\bbL^{(l)}(g_k,W_k)\;.
\end{equation}
Indeed, let us prove, for example, the (less obvious) (\ref{dEL left aux3}).
We have:
\begin{eqnarray*}
\langle d\,'_{\!1}\bbL(g_k,g_{k+1}),\eta\rangle & = & 
\left.\frac{d}{d\epsilon}\,
\bbL(g_ke^{\epsilon\eta},g_{k+1})\right|_{\epsilon=0}=
\left.\frac{d}{d\epsilon}\,
\bbL^{(l)}(g_ke^{\epsilon\eta},e^{-\epsilon\eta}W_k)\right|_{\epsilon=0}
\\ \\ & = & \langle d\,'_{\!g}\bbL^{(l)}(g_k,W_k),\eta\rangle-
\langle d_W\bbL^{(l)}(g_k,W_k),\eta\rangle\;.
\end{eqnarray*}
It remains to substitute (\ref{dEL left aux2}), (\ref{dEL left aux3}) into 
(\ref{dEL left aux1}). Taking into account that
\[
d_W\bbL^{(l)}(g_k,W_k)={\rm Ad}^*\,W_k^{-1}\cdot
d\,'_{\!W}\bbL^{(l)}(g_k,W_k)\;,
\]
we find (\ref{dEL left Ham}). Finally, notice that the notation 
(\ref{d M thru L}) is consistent with the definitions (\ref{dPi}), (\ref{d M}). 
Indeed, from these definitions it follows: $M_k=d\,'_{\!2}\bbL(g_{k-1},g_k)$, 
and the reference to (\ref{dEL left aux2}) proves (\ref{d M thru L}). 
The bracket (\ref{PB left}) is the pull--back of the standard symplectic 
bracket on $T^*G$ under the trivialization map (\ref{d M}). \qed
\vspace{2mm}

{\bf The second proof.} This proof will show us that (\ref{dEL left Ham})
describe extremals $(g_k,W_k)$ of a constrained variational principle, 
with the functional (\ref{funct left triv}), with the admissible variations
given in (\ref{variation left}). Indeed, introducing the small parameter
$\epsilon$ explicitly, i.e. replacing $\eta_k$ by $\epsilon\eta_k$, and 
writing $(g_k(\epsilon),W_k(\epsilon))$ for
$(\tilde{g}_k,\tilde{W}_k)$, we are looking for the extremum of the
functional
\[
\bbS^{(l)}(\epsilon)=\sum_{k_0}^{k_1-1}\bbL^{(l)}(g_k(\epsilon),W_k(\epsilon))\;.
\]
Considering the necessary condition
$d\bbS^{(l)}(\epsilon)/d\epsilon|_{\epsilon=0}=0$, we find:
\begin{eqnarray*}
0 & = & \sum_k\Big\langle d\,'_{\!g}\bbL^{(l)}(g_k,W_k)\,,\eta_k\Big\rangle +
\sum_k\Big\langle d\,'_{\!W}\bbL^{(l)}(g_k,W_k)\,,\eta_{k+1}-{\rm Ad}\,W_k^{-1}
\cdot\eta_k\Big\rangle\\ \\
  & = & \sum_k\Big\langle d\,'_{\!g}\bbL^{(l)}(g_k,W_k)+
d\,'_{\!W}\bbL^{(l)}(g_{k-1},W_{k-1})-{\rm Ad^*}\,W_k^{-1}\cdot
d\,'_{\!W}\bbL^{(l)}(g_k,W_k),\eta_k\Big\rangle\;.
\end{eqnarray*}
A reference to the arbitraryness of the sequence $\{\eta_k\}$ finishes the proof. 
\qed

\subsection{Reduction of left invariant Lagrangians}

Let us describe the context leading to the (discrete) Euler--Poincar\'e
equations.

Let $\Phi: G\times V\mapsto V$ be a representation of
a Lie group $G$ in a linear space $V$; we denote it by
\[
\Phi(g)\cdot v \quad\text{for} \quad g\in G\;,\;\; v\in V\;.
\]
We denote also by $\phi$ the corresponding representation of the Lie algebra 
$\gog$ in $V$:
\begin{equation}
\phi(\xi)\cdot v=\left.\frac{d}{d\epsilon}\Big(\Phi(e^{\epsilon\xi})\cdot
v\Big)\right|_{\epsilon=0} \quad \text{for} \quad \xi\in\gog\;, \;\; v\in V\;.
\end{equation}
The map $\phi^*:\gog\times V^*\mapsto V^*$  defined by
\begin{equation}\label{phi star}
\langle \phi^*(\xi)\cdot y,v\rangle=\langle y,\phi(\xi)\cdot v\rangle\qquad
\forall v\in V\;,\; y\in V^*\;,\;\xi\in\gog\;,
\end{equation}
is an anti--representation of the Lie algebra $\gog$ in $V^*$.
We shall use also the bilinear operation 
$\diamond\,: V^*\times V\mapsto \gog^*$ introduced in [HMR, CHMR] and defined as
follows: let $v\in V$, $y\in V^*$, then 
\begin{equation}\label{diamond op}
\langle y\diamond v, \xi\rangle=-\langle y,\phi(\xi)\cdot v\rangle\qquad
\forall \xi\in\gog\;.
\end{equation}
(Notice that the pairings on the left--hand side and on the right--hand side 
of the latter equation are defined on different spaces).
\vspace{2mm}

Fix an element $a\in V$, and consider the isotropy subgroup $G^{[a]}$ of $a$, 
i.e.
\begin{equation}\label{G[a]}
G^{[a]}=\{h\,:\;\Phi(h)\cdot a=a\}\subset G \;.
\end{equation}  
Suppose that the Lagrange function $\bbL(g,\hat{g})$ is invariant under 
the action of $G^{[a]}$ on $G\times G$ induced by {\it left} translations 
on $G$:
\begin{equation}\label{d left action}
\bbL(hg,h\hat{g})=\bbL(g,\hat{g})\;, \quad h\in G^{[a]}\;.
\end{equation}
The corresponding invariance property of $\bbL^{(l)}(g,W)$ is expressed as:
\begin{equation}\label{d left action for red}
\bbL^{(l)}(hg,W)=\bbL^{(l)}(g,W)\;, \quad h\in G^{[a]}\;.
\end{equation}
We want to reduce the Euler--Lagrange equations with respect to this left 
action. As a section $(G\times G)/G^{[a]}$ we choose the set $G\times O_{a}$, 
where $O_{a}$ is the orbit of $a$ under the action $\Phi$:
\begin{equation}\label{orbit}
O_a=\{\Phi(g)\cdot a\,,\;g\in G\}\subset V\;.
\end{equation} 
The reduction map is
\begin{equation}
(g,W)\in G\times G\;\mapsto\; (W,P)\in G\times O_a\;, \qquad \text{where}
\qquad P=\Phi(g^{-1})\cdot a\;,
\end{equation}
so that the reduced Lagrange function $\Lambda^{(l)}\,:\,G\times O_a
\mapsto{\mathbb R}$ is defined as 
\begin{equation}\label{left red Lagr}
\Lambda^{(l)}(W,P)=\bbL^{(l)}(g,W)\;,\quad{\rm where}\quad
 P=\Phi(g^{-1})\cdot a\;.
\end{equation}
The reduced Lagrangian $\Lambda^{(l)}(W,P)$ is well defined, because from
\[
P=\Phi(g_1^{-1})\cdot a=\Phi(g_2^{-1})\cdot a
\]
there follows $\Phi(g_2g_1^{-1})\cdot a=a$, so that $g_2g_1^{-1}\in G^{[a]}$, 
and $\bbL^{(l)}(g_1,W)=\bbL^{(l)}(g_2,W)$.
 
\begin{theorem}
{\rm a)} Consider the reduction $(g,W)\mapsto(W,P)$. The reduced Euler--Lagrange 
equations {\rm (\ref{dEL left Ham})} become the following {\bf discrete 
Euler--Poincar\'e equations}:
\begin{equation}\label{dEL left Ham red}
\left\{\begin{array}{l}
{\rm Ad}^*\,W_k^{-1}\cdot M_{k+1}=M_k+\nabla_P
\Lambda^{(l)}(W_k,P_k)\diamond P_k\;,\\ \\
P_{k+1}=\Phi(W_k^{-1})\cdot P_k\;,\end{array}\right.
\end{equation}
where
\begin{equation}\label{d M thru L red}
M_k=d\,'_{\!W}\Lambda^{(l)}(W_{k-1},P_{k-1})\in\gog^*\;.
\end{equation}
They describe extremals of the constrained variational principle, with
the functional
\begin{equation}\label{funct left triv red}
S^{(l)}=\sum_{k_0}^{k_1-1}\Lambda^{(l)}(W_k,P_k)\;,
\end{equation}
and the admissible variations $\{(\widetilde{W}_k,\widetilde{P}_k)\}$ of 
$\{(W_k,P_k)\}$ of the form
\begin{equation}\label{variation left red}
\widetilde{W}_k=W_k e^{\eta_{k+1}-{\rm Ad}\, W_k^{-1}\cdot \eta_k}\;, \qquad 
\widetilde{P}_k=P_k-\phi(\eta_k)\cdot P_k\;,
\end{equation}
where $\{\eta_k\}_{k=k_0}^{k_1}$ is an arbitrary sequence of elements of the
Lie algebra $\gog$ with $\eta_{k_0}=\eta_{k_1}=0$.

{\rm b)} If the ``Legendre transformation''
\begin{equation}\label{dLegendre left red}
(W_{k-1},P_{k-1})\in G\times O_a\mapsto 
(M_k,P_k)\in \gog^*\times O_a\;,
\end{equation}
where $P_k=\Phi(W_{k-1}^{-1})\cdot P_{k-1}$, is invertible, then
{\rm(\ref{dEL left Ham red})} define a map $(M_k,P_k)\mapsto(M_{k+1},P_{k+1})$ 
of  $\gog^*\times O_a$ which is Poisson with respect to the Poisson 
bracket 
\begin{equation}\label{PB left red}
\{F_1,F_2\}=\langle M,[\nabla_M F_1,\nabla_M F_2]\,\rangle+
\langle \nabla_P F_1,\phi(\nabla_M F_2)\cdot P\,\rangle-
\langle \nabla_P F_2,\phi(\nabla_M F_1)\cdot P\,\rangle
\end{equation}
for two arbitrary functions 
$F_{1,2}(M,P):\gog^*\times O_a\mapsto{\mathbb R}$. 
\end{theorem}
{\bf Proof} is a consequence of the following formula: if $f:G\mapsto
{\mathbb R}$ is a pull--back of the function $F:O_a\mapsto {\mathbb R}$, i.e.
\[
f(g)=F(P)=F(\Phi(g^{-1})\cdot a)\;,
\]
then
\begin{equation}\label{aux1}
d\,'f(g)=\nabla_P F(P)\diamond P\;.
\end{equation}
(Indeed,
\[
\langle d\,'f(g),\xi\rangle=\left.\frac{d}{d\epsilon}f(ge^{\epsilon\xi})\right|
_{\epsilon=0}=\left.\frac{d}{d\epsilon}F(\Phi(e^{-\epsilon\xi})\cdot P)\right|
_{\epsilon=0}=-\langle \nabla_PF(P),\phi(\xi)\cdot P\rangle=
\langle \nabla_PF(P)\diamond P,\xi\rangle\;;
\]
the last equality is the definition (\ref{diamond op})). In particular,
plugging
\[
d\,'_{\!g}\bbL^{(l)}=\nabla_P\Lambda^{(l)}\diamond P\;,\qquad
d\,'_{\!W}\bbL^{(l)}=d\,'_{\!W}\Lambda^{(l)}
\]
into the first equation in (\ref{dEL left Ham}), we come to the first equation
in (\ref{dEL left Ham red}). Similarly, the Poisson bracket (\ref{PB left}) 
turns into (\ref{PB left red}) by use of (\ref{aux1}). \qed
\vspace{2mm}

{\bf Remark 1.} The formula (\ref{PB left red}) defines a Poisson bracket 
not only on $\gog^*\times O_a$, but on all of $\gog^*\times V$. Rewriting
this formula as 
\begin{equation}\label{PB left red once more}
\{F_1,F_2\}=\langle M,[\nabla_M F_1,\nabla_M F_2]\,\rangle+
\langle  P, \phi^*(\nabla_M F_2)\cdot \nabla_P F_1-
\phi^*(\nabla_M F_1)\cdot \nabla_P F_2\,\rangle
\end{equation}
one immediately identifies this bracket with the Lie--Poisson bracket of the 
semiproduct Lie algebra $\gog\ltimes V^*$ corresponding to the representation
$-\phi^*$ of $\gog$ in $V^*$.
\vspace{2mm}

{\bf Remark 2.} In an important particular case of constructions of this 
section, the vector space is chosen as the Lie algebra 
of our basic Lie group: $V=\gog$, the group representation is the adjoint one:
$\Phi(g)\cdot v={\rm Ad\; g}\cdot v$, so that $\phi(\xi)\cdot v={\rm ad\;}\xi
\cdot v=[\xi,v]$, and the bilinear operation $\diamond$ is nothing but
the coadjoint action of $\gog$ on $\gog^*$: 
$y\diamond v={\rm ad}^*\, v\cdot y$. This is the framework, e.g., for the
heavy top mechanics. The model studied in [BS] belongs to this class.

\setcounter{equation}{0}
\section{Right symmetry reduction}
All constructions here are parallel to those of the previous section,
so we restrict ourselves to formulations of the basic results only.

\subsection{Right trivialization}

Consider the right trivialization of the cotangent bundle $T^*G$:
\begin{equation}\label{d right triv *}
(g_k,m_k)\in G\times\gog^*\;\mapsto \;(g_k,\Pi_k)\in T^*G\;,
\end{equation}
where 
\begin{equation}\label{d m}
\Pi_k=R_{g_k^{-1}}^* m_k \quad\Leftrightarrow\quad m_k=R_{g_k}^*\Pi_k\;.
\end{equation}
Consider also the map
\begin{equation}\label{d right triv}
(g_k,w_k)\in G\times G\; \mapsto \;(g_k,g_{k+1})\in G\times G\;,
\end{equation}
where
\begin{equation}\label{w}
g_{k+1}=w_kg_k \quad\Leftrightarrow\quad w_k=g_{k+1}g_k^{-1}\;.
\end{equation}
Denote the pull--back of the Lagrange function under (\ref{d right triv}) 
through
\begin{equation}\label{dLagr right}
\bbL^{(r)}(g_k,w_k)=\bbL(g_k,g_{k+1})\;.
\end{equation}
\begin{proposition}
Consider the functional 
\[
\bbS^{(r)}=\sum_{k_0}^{k_1-1}\bbL^{(r)}(g_k,w_k)\;,
\] 
and its extremals with respect to variations $\{(\tilde{g}_k,\tilde{w}_k)\}$
of $\{(g_k,w_k)\}$ of the form
\[
\tilde{g}_=e^{\eta_k}g_k\;, \qquad
\tilde{w}_k=e^{\eta_{k+1}-{\rm Ad}\,w_k\cdot\eta_k}
 w_k\;,
\]
where $\{\eta_k\}_{k=k_0}^{k_1}$ is an arbitrary sequence of elements of $\gog$
with $\eta_{k_0}=\eta_{k_1}=0$. The difference equations for extremals of 
this constrained variational problem read:
\begin{equation}\label{dEL right Ham}
\left\{\begin{array}{l}
{\rm Ad}^*\,w_k\cdot m_{k+1}=m_k+d_g\bbL^{(r)}(g_k,w_k)\;,\\ \\
g_{k+1}=w_kg_k\;,\end{array}\right.
\end{equation}
where
\begin{equation}\label{d m thru L}
m_k=d_w\bbL^{(r)}(g_{k-1},w_{k-1})\in\gog^*\;.
\end{equation}
If the ``Legendre transformation'' 
\begin{equation}\label{dLegendre right}
(g_{k-1},w_{k-1})\in G\times G\mapsto (g_k,m_k)\in G\times\gog^*\;,
\end{equation}
where $g_k=w_{k-1}g_{k-1}$, is invertible, then {\rm(\ref{dEL right Ham})} 
define a map $(g_k,m_k)\mapsto(g_{k+1},m_{k+1})$ which is symplectic with 
respect to the following Poisson bracket on $G\times\gog^*$:
\begin{equation}\label{PB right}
\{f_1,f_2\}=-\langle d_g f_1,\nabla_m f_2\rangle+
\langle d_g f_2,\nabla_m f_1\rangle-
\langle m,[\nabla_m f_1,\nabla_m f_2]\,\rangle\;.
\end{equation}
\end{proposition}
{\bf Proof.} This time the discrete Euler--Lagrange equations (\ref{dEL})
are rewritten  as
\begin{equation}\label{dEL right aux1}
d_1\bbL(g_k,g_{k+1})+d_2\bbL(g_{k-1},g_k)=0\;,
\end{equation}
and the expressions for these Lie derivatives in terms of $(g,w)$ read: 
\begin{equation}\label{dEL right aux2}
d_2\bbL(g_{k-1},g_k)=d_w\bbL^{(r)}(g_{k-1},w_{k-1})\;,
\end{equation}
\begin{equation}\label{dEL right aux3}
d_1\bbL(g_k,g_{k+1})=d_g\bbL^{(r)}(g_k,w_k)-
d\,'_{\!w}\bbL^{(r)}(g_k,w_k)=d_g\bbL^{(r)}(g_k,w_k)-
{\rm Ad}^*\,w_k\cdot d_w\bbL^{(r)}(g_k,w_k)\;.
\end{equation}
Finally, the expression (\ref{d m thru L}) is consistent with the definitions
(\ref{dPi}), (\ref{d m}), which imply that $m_k=d_2\bbL(g_{k-1},g_k)$,
and a reference to (\ref{dEL right aux2}) finishes the proof. \qed

\subsection{Reduction of right invariant Lagrangians}
 Assume that the function 
$\bbL^{(r)}$ is invariant under the action of $G^{[a]}$ on 
$G\times G$ induced by {\it right} translations on $G$:
\begin{equation}\label{d right action for red}
\bbL^{(r)}(gh,w)=\bbL^{(r)}(g,w)\;, \quad h\in G^{[a]}\;.
\end{equation}
Define the reduced Lagrange function $\Lambda^{(r)}\,:\,G\times O_a
\mapsto{\mathbb R}$ as
\begin{equation}\label{right red dLagr}
\Lambda^{(r)}(w,p)=\bbL^{(r)}(g,w)\;,\quad{\rm where}\quad
 p=\Phi(g)\cdot a\;.
\end{equation}
\begin{theorem}
{\rm a)} Consider the reduction $(g,w)\mapsto(w,p)$. The reduced 
Euler--Lagrange equations {\rm (\ref{dEL right Ham})} become the following
{\bf discrete Euler--Poincar\'e equations}:
\begin{equation}\label{dEL right Ham red}
\left\{\begin{array}{l}
{\rm Ad}^*\,w_k\cdot m_{k+1}=m_k-\nabla_p\Lambda^{(r)}(w_k,p_k)\diamond
p_k\;,\\ \\
p_{k+1}=\Phi(w_k)\cdot p_k\;,\end{array}\right.
\end{equation}
where
\begin{equation}\label{d m thru L red}
m_k=d_w\Lambda^{(r)}(w_{k-1},p_{k-1})\in\gog^*\;.
\end{equation}
They describe extremals of the constrained variational principle, with
the functional
\begin{equation}\label{funct right triv red}
S^{(r)}=\sum_{k_0}^{k_1-1}\Lambda^{(r)}(w_k,p_k)\;,
\end{equation}
and the admissible variations $\{(\widetilde{w}_k,\widetilde{p}_k)\}$ of 
$\{(w_k,p_k)\}$ of the form
\begin{equation}\label{variation right red}
\widetilde{w}_k=e^{\eta_{k+1}-{\rm Ad}\, w_k\cdot \eta_k}w_k \;, \qquad 
\widetilde{p}_k=p_k+\phi(\eta_k)\cdot p_k\;,
\end{equation}
where $\{\eta_k\}_{k=k_0}^{k_1}$ is an arbitrary sequence of elements of the
Lie algebra $\gog$ with $\eta_{k_0}=\eta_{k_1}=0$.

{\rm b)} If the ``Legendre transformation'' 
\begin{equation}\label{dLegendre right red}
(w_{k-1},p_{k-1})\in G\times O_a\mapsto 
(m_k,p_k)\in \gog^*\times O_a\;,
\end{equation}
where $p_k=\Phi(w_{k-1})\cdot p_{k-1}$, is invertible, then {\rm(\ref
{dEL right Ham red})} define a map $(m_k,p_k)\mapsto(m_{k+1},p_{k+1})$ of  
$\gog^*\times O_a$ which is Poisson with respect to the bracket 
\begin{equation}\label{PB right red}
\{F_1,F_2\}=-\langle m,[\nabla_m F_1,\nabla_m F_2]\,\rangle
-\langle \nabla_p F_1,\phi(\nabla_m F_2)\cdot p\,\rangle+
\langle \nabla_p F_2,\phi(\nabla_m F_1)\cdot p\,\rangle
\end{equation}
for two arbitrary functions 
$F_{1,2}(m,p):\gog^*\times O_a\mapsto{\mathbb R}$. This formula indeed
defines a Poisson bracket on all of $\gog^*\times V$, the Lie--Poisson
bracket of the semiproduct Lie algebra $\gog\ltimes V^*$ corresponding
to the representation $-\phi^*$ of $\gog$ in $V^*$.
\end{theorem}
{\bf Proof} is based on the following simple result: if
$f(g)=F(\Phi(g)\cdot a)$, then
\[
df(g)=- \nabla_p F(p)\diamond p\;.
\]
Notice that the brackets (\ref{PB left red}) and (\ref{PB right red})
essentially coincide (differ only by a sign). \qed
\vspace{2mm}

A table summarizing the unreduced and reduced Lagrangian equations of motion, 
both in the continuous and discrete time formulations, is put in Appendix A.
The continuous time results were obtained in [HMR],[CHMR]; they also may
be derived by taking a continuous limit of our present results.

\section{Conclusion}
We consider the discrete time Lagrangian mechanics on Lie groups
as an important source of symplectic and, more general, Poisson maps.
Moreover, from some points of view the variational (Lagrangian) structure
is even more fundamental and important than the Poisson (Hamiltonian) one
(cf. \cite{HMR}, \cite{MPS}, where a similar viewpoint is represented).
In particular, discrete Lagrangians on $G\times G$ may serve as models
for the rigid body motion (cf. \cite{WM}). The integrable cases of Euler (a
free rotation of a rigid body fixed at the center of mass) and of
Lagrange (symmetric spinning top) are discretized in this framework 
preserving the integrability property in [V],[MV] and in [BS], respectively. 
It would be interesting and important to apply the above theory to the 
infinite dimensional case, e.g. to discretization of ideal compressible 
fluids motion (see [HMR]).

\begin{appendix}

\setcounter{equation}{0}
\section{Euler--Lagrange and Euler--Poincar\'e equations}

\noindent {\small
\begin{tabular}{|l|l|}\hline
 & \\
\multicolumn{1}{|c|}{\qquad CONTINUOUS TIME\qquad\;} & 
\multicolumn{1}{|c|}{\quad\qquad\qquad\qquad DISCRETE TIME
\qquad\qquad\qquad\qquad}\\ 
 & \\ \hline\hline   
\multicolumn{2}{|c|}{}\\ 
\multicolumn{2}{|c|}{General Lagrangian systems}\\ 
\multicolumn{2}{|c|}{} \\ \hline
 & \\
$\bL(g,\dot{g})$ & $\bbL(g_k,g_{k+1})$ \\ 
 & \\
$\left\{\begin{array}{l} \Pi=\nabla_{\dot{g}}\bL\\
                                           \dot{\Pi}=\nabla_g\bL
                          \end{array}\right.$ &
                $\left\{\begin{array}{l} \Pi_k=-\nabla_1\bbL(g_k,g_{k+1})\\
                                        \Pi_{k+1}=\nabla_2\bbL(g_k,g_{k+1})
                          \end{array}\right.$ \\
 & \\
\hline\hline 
\end{tabular}

\noindent
\begin{tabular}{|l|l|}\hline\hline
\multicolumn{2}{|c|}{} \\
\multicolumn{2}{|c|}{Left trivialization, left symmetry reduction: 
$\;\;M=L_g^*\Pi\,,\;\; P=\Phi(g^{-1})\cdot a$} \\
\multicolumn{2}{|c|}{} \\ \hline
 & \\
$\bL(g,\dot{g})=\cL^{(l)}(\Omega,P)$ &
$\bbL(g_k,g_{k+1})=\Lambda^{(l)}(W_k,P_k)$ \\ 
$\Omega=L_{g^{-1}*}\dot{g}\,,\;\;P=\Phi(g^{-1})\cdot a$ & 
$W_k=g_k^{-1}g_{k+1}\,,\;\;P_k=\Phi(g_k^{-1})\cdot a$\\ 
 & \\
$M=L_g^*\Pi=\nabla_{\Omega}\cL^{(l)}$ & 
$M_k=L_{g_k}^*\Pi_k=d\,'_{\!W}\Lambda^{(l)}(W_{k-1},P_{k-1})$ \\ 
 & \\
$\left\{\begin{array}{l}
\dot{M}={\rm ad}^*\,\Omega\cdot M+\nabla_P\cL^{(l)}\diamond P\\ 
\dot{P}=-\phi(\Omega)\cdot P\end{array}\right.$ &
 $\left\{\begin{array}{l}
{\rm Ad}^*\,W_k^{-1}\cdot M_{k+1}=M_k+
\nabla_P\Lambda^{(l)}(W_k,P_k)\diamond P_k\\ 
P_{k+1}=\Phi(W_k^{-1})\cdot P_k\end{array}\right.$ \\
 & \\
\hline\hline 
\multicolumn{2}{|c|}{} \\
\multicolumn{2}{|c|}{Right trivialization, right symmetry reduction:
$\;\;m=R_g^*\Pi\,,\;\; p=\Phi(g)\cdot a$} \\ 
\multicolumn{2}{|c|}{} \\ \hline
 & \\
$\bL(g,\dot{g})=\cL^{(r)}(\omega,p)$ &
$\bbL(g_k,g_{k+1})=\Lambda^{(r)}(w_k,p_k)$ \\ 
$\omega=R_{g^{-1}*}\dot{g}\,,\;\;p=\Phi(g)\cdot a$ & 
$w_k=g_{k+1}g_k^{-1}\,,\;\;p_k=\Phi(g_k)\cdot a$\\ 
 & \\
$m=R_g^*\Pi=\nabla_{\omega}\cL^{(r)}$ & 
$m_k=R_{g_k}^*\Pi_k=d_w\Lambda^{(r)}(w_{k-1},p_{k-1})$ \\ 
 & \\
$\left\{\begin{array}{l}
\dot{m}=-{\rm ad}^*\,\omega\cdot m-\nabla_p\cL^{(r)}\diamond p\\ 
\dot{p}=\phi(\omega)\cdot p\end{array}\right.$ &
 $\left\{\begin{array}{l}
{\rm Ad}^*\,w_k\cdot m_{k+1}=m_k-
\nabla_p\Lambda^{(r)}(w_k,p_k)\diamond p_k\\ 
p_{k+1}=\Phi(w_k)\cdot p_k\end{array}\right.$ \\
 & \\
\hline\hline
\end{tabular}
}
\vspace{5mm}

The relation between the continuous time and the discrete time equations 
is established, if we set
\[
g_k=g\;,\qquad g_{k+1}=g+\varepsilon\dot{g}+O(\varepsilon^2)\;,\qquad
\bbL(g_k,g_{k+1})=\varepsilon\bL(g,\dot{g})+O(\varepsilon^2)\;;
\]
\[
P_k=P\;,\qquad W_k=\ed+\varepsilon\Omega+O(\varepsilon^2)\;,\qquad
\Lambda^{(l)}(W_k,P_k)=\varepsilon\cL^{(l)}(\Omega,P)+O(\varepsilon^2)\;;
\]
\[
p_k=p\;,\qquad w_k=\ed+\varepsilon\omega+O(\varepsilon^2)\;,\qquad
\Lambda^{(r)}(w_k,p_k)=\varepsilon\cL^{(r)}(\omega,p)+O(\varepsilon^2)\;.
\]

\setcounter{equation}{0}
\section{Notations}
We fix here some notations and definitions used throughout the paper.

Let $G$ be a Lie group with the Lie algebra $\gog$, and let $\gog^*$ be a dual 
vector space to $\gog$. We identify $\gog$ and $\gog^*$ with the tangent space 
and the cotangent space to $G$ in the group unit, respectively:
\[
\gog=T_e G\;,\qquad \gog^*=T^*_e G\;.
\]
The pairing between the cotangent and the tangent spaces $T^*_g G$ and 
$T_g G$ in an arbitrary point $g\in G$ is denoted by 
$\langle\cdot,\cdot\rangle$. The left and right translations in the group
are the maps $L_g\,,R_g\,:G\mapsto G$ defined by
\[
L_g\,h=gh\;,\qquad R_g\,h=hg\qquad\forall h\in G\;,
\]
and $L_{g*}\,$, $R_{g*}$ stand for the differentials of these maps:
\[
L_{g*}\,:\,T_h G\mapsto T_{gh} G\;, \qquad R_{g*}\,:\,T_h G\mapsto T_{hg} G\;.
\]
We denote by
\[
{\rm Ad}\,g = L_{g*}R_{g^{-1}*}\,:\gog\mapsto\gog
\]
the adjoint action of the Lie group $G$ on its Lie algebra $\gog=T_e G$. 
The linear operators
\[
L_g^*\,:\, T^*_{gh} G\mapsto T^*_h G\;, \qquad
R_g^*\,:\, T^*_{hg} G\mapsto T^*_h G
\]
are conjugated to $L_{g*}\,$, $R_{g*}\,$, respectively, via the pairing
$\langle\cdot,\cdot\rangle$:
\[
\langle L_g^*\xi,\eta\rangle=\langle\xi,L_{g*}\eta\rangle\quad{\rm for}
\quad \xi\in T^*_{gh} G\,,\;\; \eta\in T_h G\;,
\]
\[
\langle R_g^*\xi,\eta\rangle=\langle\xi,R_{g*}\eta\rangle\quad{\rm for}
\quad \xi\in T^*_{hg} G\,,\;\; \eta\in T_h G\;.
\]
The coadjoint action of the group
\[
{\rm Ad}^*\,g = L^*_{g}R^*_{g^{-1}}\,:\gog^*\mapsto\gog^*
\]
is conjugated to ${\rm Ad}\,g$ via the pairing $\langle\cdot,\cdot\rangle$:
\[
\langle {\rm Ad}^*\, g\cdot\xi,\eta\rangle=\langle\xi, 
{\rm Ad}\,g\cdot\eta\rangle\quad{\rm for}
\quad \xi\in\gog^*\,,\;\; \eta\in\gog\;.
\]
The differentials of ${\rm Ad}\,g$ and of ${\rm Ad}^*\,g$ with respect to $g$
in the group unity $e$ are the operators 
\[
{\rm ad}\,\eta\,:\gog\mapsto\gog \qquad{\rm and}\qquad
{\rm ad}^*\,\eta\,:\gog^*\mapsto\gog^*\;,
\]
respectively, also conjugated via the pairing $\langle\cdot,\cdot\rangle$:
\[
\langle{\rm ad}^*\,\eta\cdot\xi,\zeta\rangle=\langle\xi,{\rm ad}\,\eta\cdot\zeta
\rangle \qquad \forall \xi\in\gog^*\,,\;\;\zeta\in\gog\;.
\]
The action of ad is given by applying the Lie bracket in $\gog$:
\[
{\rm ad}\,\eta\cdot\zeta=[\eta,\zeta]\;, \quad \forall\zeta\in\gog\;.
\]

Finally, we need the notion of gradients of functions on vector spaces
and on manifolds. If $\cX$ is a vector space, and $f:\cX\mapsto{\mathbb R}$ is a 
smooth function, then the gradient $\nabla f:\cX\mapsto\cX^*$ is defined via
the formula
\[
\langle\nabla f(x),y\rangle=\left.\frac{d}{d\epsilon}\,f(x+\epsilon y)
\right|_{\epsilon=0}\;, \qquad \forall y\in\cX\;.
\]
Similarly, for a function $f:G\mapsto{\mathbb R}$ on a smooth manifold $G$ its
gradient $\nabla f\,:\,G\mapsto T^* G$ is defined in the following way: for an
arbitrary $\dot{g}\in T_g G$ let $g(\epsilon)$ be a curve in $G$ through
$g(0)=g$ with the tangent vector $\dot{g}(0)=\dot{g}$. Then
\[
\langle\nabla f(g),\dot{g}\rangle=\left.\frac{d}{d\epsilon}\,f(g(\epsilon))
\right|_{\epsilon=0}\;.
\]
If $G$ is a Lie group, then two convenient ways to define a curve in $G$ through
$g$ with the tangent vector $\dot{g}$ are the following:
\[
g(\epsilon)=e^{\epsilon\eta}g\;,\quad \eta=R_{g^{-1}*}\,\dot{g}\;,
\]
and
\[
g(\epsilon)=ge^{\epsilon\eta}\;,\quad \eta=L_{g^{-1}*}\,\dot{g}\;,
\]
which allows to establish the connection of the gradient $\nabla f$ with
the (somewhat more convenient) notions of the left and the right Lie derivatives
of a function $f:G\mapsto{\mathbb R}$:
\[
\nabla f(g)=R_{g^{-1}}^*\,df(g)=L_{g^{-1}}^*\,d\,'f(g)\;.
\] 
Here $df:G\mapsto\gog^*$ and $d\,'f:G\mapsto\gog^*$ are defined via 
the formulas
\[
\langle df(g),\eta\rangle=\left.\frac{d}{d\epsilon}\,f(e^{\epsilon\eta}g)
\right|_{\epsilon=0}\;,\qquad \forall \eta\in\gog\;,
\]
\[
\langle d\,'f(g),\eta\rangle=\left.\frac{d}{d\epsilon}\,f(ge^{\epsilon\eta})
\right|_{\epsilon=0}\;,\qquad \forall \eta\in\gog\;.
\]

\end{appendix}

\end{document}